\documentclass[12pt]{amsart}
\usepackage{epsfig}
\usepackage{amscd}

\newtheorem{theorem}{Theorem}
\newtheorem{proposition}{Proposition}

\newtheorem{fact}{Fact}

\newcommand{\rmk}{\hfill\mbox{}\par}
\newcommand{\lmk}{\noindent\mbox{}\hfill}
\newcommand{\cL}{{\mathcal L}}
\newcommand{\bC}{{\hspace{-0.3pt}\mathbb C}\hspace{0.4pt}}
\newcommand{\lcr}{\raisebox{-5pt}{\mbox{}\hspace{1pt}
                 \includegraphics{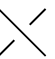}\hspace{1pt}\mbox{}}}
\newcommand{\ift}{\raisebox{-5pt}{\mbox{}\hspace{1pt}
                 \includegraphics{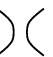}\hspace{1pt}\mbox{}}}
\newcommand{\zer}{\raisebox{-5pt}{\mbox{}\hspace{1pt}
                 \includegraphics{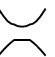}\hspace{1pt}\mbox{}}}
\newcommand{\kt}{K_t}
\newcommand{\g}{SL_2(\bC)} 
\newcommand{\un}{\underline}
\newcommand{\cA}{{\mathcal A}}
\newcommand{\qg}{{_qSL_2}}
\newcommand{\ov}{\overline}
\newcommand{\mat}{\ov{\cA}_t}
\newcommand{\bA}{{\hspace{-0.3pt}\mathbb A}\hspace{0.3pt}}
\newcommand{\cO}{{\mathcal O}}
\newcommand{\cG}{{\mathcal G}}

\title{The Kauffman Bracket Skein as an Algebra of Observables}

\author{DOUG BULLOCK}

\address{Department of Mathematics, Boise State University,
Boise, ID 83725}

\email{bullock@math.idbsu.edu}

\author{CHARLES FROHMAN}

\address{Department of Mathematics, University of Iowa, Iowa City, IA
52242}

\email{frohman@math.uiowa.edu}

\author{JOANNA KANIA-BARTOSZY\'{N}SKA}

\address{Department of Mathematics, Boise State University, Boise, ID
83725}

\email{kania@math.idbsu.edu}

\thanks{This research was partially supported by an NSF-DMS Postdoctoral
Research Fellowship, and by NSF grants DMS-9803233
 and DMS-9971905.}
  
\begin{document}

\begin{abstract} 
We prove that the Kauffman bracket skein algebra of a cylinder over
a surface with boundary, defined over complex
numbers, is isomorphic to the observables of an appropriate lattice
gauge field theory. 

\end{abstract}

\maketitle

\section{Introduction}

Lattice gauge field theory brings the representation theory of an
underlying manifold and its quantum invariants into the same setting. 
Consider the case of a cylinder over a compact, oriented surface with
boundary. A 
lattice model of the surface determines an algebra of gauge invariant
fields (i.e. observables). In the classical case, based on a connected,
simply connected Lie group $G$, observables are the
characters of the fundamental group of the lattice represented
in $G$. Wilson loops can be understood as traces of conjugacy classes
in the fundamental group of the lattice. 
For 
the theory based on a Drinfeld-Jimbo deformation of a simple Lie
algebra $\mathfrak{g}$, the observables are a deformation quantization
of the $G$-characters of the surface with respect to the standard
Poisson 
structure \cite{lattice}. In the case of $U_h(sl_2)$ this, together
with the classical isomorphism  \cite{bul,przysik},  allowed us to
prove that  the algebra of observables is the Kauffman
bracket skein algebra of the surface, completed as an algebra  over
formal power series.

In this paper we return to an analytic setting in which the
deformation parameter is any complex number other than a root of unity.
The analogous theorem relating observables and the Kauffman
bracket skein algebra is again true. The proof is based on the
combinatorial equivalence between the Temperley-Lieb algebra
and the quantized invariant theory of $SL_2$; it
does not explicitly use the relationship with surface characters.

The paper is organized as follows. Section \ref{kauffman} recalls
definitions and associated formulas of the Kauffman bracket skein
algebra.  Section \ref{reps} summarizes, for the generic parameter,
the construction of a quasi-triangular matrix model of quantum $SL_2$.
Section \ref{lgft} outlines basic definitions and constructions of
lattice gauge field theory. Section \ref{temperley} describes the
correspondence between skeins and intertwiners in the Verlinde algebra
for quantum $SL_2$.
Finally, Section \ref{observable} contains a proof of the main
theorem.

\section{Kauffman Bracket Skein Algebra}\label{kauffman}

Let $t$ be a complex number that is neither $0$ nor a root of unity.
Suppose that $F$ is a compact, oriented surface with boundary and
$I$ is a closed interval. Denote by $\cL$ the set of isotopy classes 
of framed links in $F\times I$, including the empty link. 
Let $\bC\cL$ be the vector space with basis
$\cL$.   Define $S_t$ to be the
subspace of $\bC \cL$ spanned by all expressions of the form
$\displaystyle{\lcr-t\zer-t^{-1}\ift}$ and $\bigcirc+t^2+t^{-2}$,
where the framed links in each expression are identical outside the
balls pictured in the diagrams.

The Kauffman bracket skein algebra $K_t(F\times I)$ is the quotient 
$\kt(F\times I) = \bC \cL / S_t$. Multiplication is 
given by laying one link over the other. More
precisely, if $\alpha$ and $\beta$ are in $\cL$, isotop
them so that $\alpha$ lies in  $F\times [0,\frac{1}{2})$ ,
and $\beta$ in $F\times (\frac{1}{2},1]$.  
Then $\alpha*\beta$ is the union of
these two links in $F\times [0,1]$. Extend linearly to a product on
$\bC\cL$. Since $S_t$ is an ideal, the product descends, making the skein
module into a skein algebra.
Since the algebra structure depends on the specific product structure
of $F\times I$, rather than its topological type, we  use the
notation $\kt(F)$.

We use the standard convention of modeling a skein in $\kt(F)$ on a framed,
admissibly colored, trivalent graph. An admissible coloring  is an assignment of a
nonnegative integer to each edge so that the colors at each
vertex form admissible triples. A triple $(a,b,c)$ is admissible if
$a\leq b+c$, $b\leq a+c$,  
$c\leq a+b$ and $a+b+c$ is even.  The corresponding
skein in $\kt(F)$ is obtained by replacing each edge labeled with the
letter $m$ by the $m$-th Jones--Wenzl
idempotent (see \cite{We}, or \cite[p.\ 136]{Li}), and replacing
trivalent vertices with Kauffman triads (see \cite[Fig.\ 14.7]{Li}).
If $s$ is a trivalent spine of $F$ then the set of skeins carried by
admissible colorings of $s$   forms a basis ${\mathcal B}_1$ for
$\kt(F)$.
If $F$ is an annulus, ${\mathcal B}_1$ consists of skeins obtained by
labeling the core with a Jones--Wenzl idempotent. One may think of the
core as a ``trivalent'' spine with one vertex, whose admissible labels
are $\{(n,n,0)\}$.
The space $\kt(F)$ also has a basis ${\mathcal B}_2$ consisting of all
links with 
simple diagrams on $F$, i.e.\ with no crossings and no trivial components.

\section{Representations}\label{reps}

The details of the following can be found in \cite{matrix}.
Let $\cA_t$ be the unital Hopf algebra on $X$, $Y$, $K$, $K^{-1}$,
with relations: 
\begin{align*} KX=t^2XK, \quad KY=t^{-2}YK,\\
 XY-YX= \frac{K^2-K^{-2}}{t^2-t^{-2}}, \quad KK^{-1}=1. \end{align*}
Let  $\un{m}$ denote the irreducible
$(m+1)$-dimensional representation of $\cA_t$. Fixing an ordered
basis for $\un{m}$ we define
linear functionals  ${^mc}^i_j : A_t \rightarrow \bC$ to be the
coefficient in the $i$-th row and $j$-th column in the representation
$\un{m}$. The ${^mc}^i_j$ form a basis for the stable subalgebra
$\qg$ of the Hopf algebra dual $A_t^o$. (Here $q=t^4$.)
Define
\[ \mat = \prod_{m=0}^{\infty}  M_{m+1}(\bC)\]
and give it the product topology. A typical element of $\mat$ is
a sequence of arbitrarily chosen matrices in which the $i$-th term
is an $(i+1) \times (i+1)$ matrix.

Let $\rho_m:\cA_t \rightarrow M_{m+1}(\bC)$ be the homomorphism
corresponding to the representation $\un{m}$. The homomorphism
\begin{equation}\label{Theta}
\Theta : \cA_t \rightarrow \mat,
\end{equation}
given by $\Theta(Z)=(\rho_0(Z),\rho_1(Z),\rho_2(Z), \ldots)$, is
injective and its image is dense in $\mat$ (see \cite{matrix}). The
algebra $\mat$ is the completion of $\cA_t$ by equivalence classes of
Cauchy sequences in the weak topology from $\qg$. It has the structure
of a topological ribbon Hopf algebra.  The projection of $\mat$ onto
its $(m+1)$-st factor is an irreducible representation of $\mat$, also
denoted by $\un{m}$. Composing ${^mc}^i_j$ with this projection yields
a function on $\mat$, also called ${^mc}^i_j$. Thus $\qg$ is understood
to lie in $(\mat)^o$.

\section{lattice gauge field theory}\label{lgft}

In this section we recall basic definitions and  constructions of
lattice gauge field theory. Details can be found in \cite{lattice}.

Let $\Gamma$ be an oriented, ciliated graph; i.e.\ the edges are
oriented and the vertices carry a linear ordering of the adjacent
edges.
One can think of this as instructions for building a strip
and disk model of an oriented surface $F$ having $\Gamma$ as a
strong deformation retract. Vertices correspond to the disks, 
edges to  strips, and the ciliation determines how to glue the
strips to the disks. The surface $F$ is called the envelope of
$\Gamma$.
Define a space of  connections: 
\[\bA(\Gamma) = \bigotimes_{\text{edges $e$}} (\mat)_e;\]
and define an algebra of  fields:
\[C[\bA(\Gamma)] = \bigotimes_{\text{edges $e$}}  (\qg)_e.\]

Note that fields are functions on connections in the obvious way.
The connections form a coalgebra with  comultiplication as
defined in \cite{lattice}. Multiplication
of fields is the convolution product dual to comultiplication
of connections. There is an action of the  gauge algebra,
\[\cG(\Gamma)= \bigotimes_{\text{vertices $v$}} (\mat)_v,  \]
on the space of connections, and adjointly on fields. The invariant
part of the gauge fields under this action is called the 
observables, $\cO(\Gamma)$. The multiplication of fields
restricts to make  $\cO(\Gamma)$ into an algebra, which is a
deformation of the $\g$-characters of $\pi_1(F)$.

Let $V$ be a representation of the Hopf algebra $\mat$, that is $V$ is
a finite dimensional left module over $\mat$.

The dual vector space to $V$ carries two distinct $\mat$-module
structures.
When $\mat$ acts on the left, the dual module is denoted by $V^*$.
The action is:
\[ Z\cdot\phi(v)=\phi(S(Z)\cdot v),\]
for any $Z \in \mat$, $v \in V$ and $\phi$ in the dual vector space to
$V$. When $\mat$ acts on the right, denote the dual by $V'$ with
\[\phi(v)\cdot Z =\phi(Z\cdot v).\]

There is an alternate description of observables in terms of
``colorings'' of the lattice. Assigning  a representation 
$V_e$ to each positively oriented edge $e$ of the lattice $\Gamma$ determines a map 
\[\bA(\Gamma)\to\bigotimes_e(V^*_e\otimes V_e).\]
This yields,  at each vertex, a tensor product of
representations coming from the edges adjoining that vertex taken in
the order given by the ciliation. Use the representation $V_e$ for
the edges $e$ starting at a vertex and the dual $V^*_e$ for the edges
terminating there. The resulting representation at a vertex $v$ is
denoted by $V_v$.
Finally, choose
$\phi_v\in{\rm Inv}(V'_v)$ for each vertex $v$.
The element $\otimes_v(\phi_v)$ is evaluated on a connection
$\otimes_ex_e$ by mapping $\otimes_ex_e$ to 
$\otimes_e(V^*_e\otimes V_e)$  and then re-parsing to an element of
$\otimes_vV_v$. By \cite[Corollary 1]{lattice} every observable is
a linear combination of observables of this form.

Now assume that $\Gamma$ is a trivalent lattice colored by
irreducible representations. The coloring is admissible if, at each
vertex, the integers corresponding to the colorings of the incident
edges  form an admissible triple.

\begin{proposition} \label{basis}
Suppose that $\Gamma$ is an admissibly colored trivalent lattice. For
each $V_v$ there exists a non-zero dual element invariant under the
right action of $\mat$. The tensor product of these invariants over all
vertices defines an observable.

The set of such observables, one for each admissible coloring, is a
basis for $\mathcal{O}(\Gamma)$.
\end{proposition} 

\proof
Let $c=\{\un{m}_{\mbox{}\,e}\ |\ e \mbox{ is an edge of }\Gamma\}$ be an
admissible coloring.
Note that admissibility implies a $1$-dimensional invariant subspace in
each $V'_v$. Hence there is a non-zero observable  $o_c=\otimes_v\phi_v$.
Since $o_c$ is nonzero, there exists
$X_c\in\otimes_v V_v$ so that $o_c(X_c)\neq 0$. 

Let
$P:\otimes_vV_v\rightarrow\otimes_e
\left(\un{m}_{\mbox{}\,e}^*\otimes\un{m}_{\mbox{}\,e}\right)$
be the ``parsing'' map. Let $\iota_m$ be the map
\[\un{m}^*\otimes\un{m}\cong
M_{m+1}(\bC)\rightarrow\prod_nM_n(\bC)=\mat,\] 
where the inclusion is given by forming a sequence that has all zero
matrices except for the $(m+1)$-st entry corresponding to $\un{m}$.

Define $x_c$ to be the connection
$\left(\otimes_e\iota_{m_e}\right)\left(P(X_c)\right)$. Clearly
$o_c(x_c)\neq 0$. Since $\rho_m\circ\iota_n =0$ unless $m=n$, $x_c$ is
annihilated by all observables constructed from colorings different
than $c$.

From this we conclude that any set of observables constructed from
distinct colorings is independent. By \cite{lattice} they span.
\qed

There is a map
\begin{equation}\label{Phi}
\Phi : \kt(F) \to \cO(\Gamma)
\end{equation}
that assigns to each framed link a Wilson
operator. Classically, a Wilson operator is the trace of the holonomy of a
connection along some fixed loop. The construction in the quantum
setting is described in \cite{lattice}.
In a theory based on the Drinfeld-Jimbo $U_h(sl_2)$ and its Hopf
algebra dual, the observables are isomorphic to the Kauffman bracket
skein algebra of the lattice envelope.
Since $U_h(sl_2)$ is a complete topological algebra over $\bC[[h]]$,
it was necessary to complete the skein algebra as well.
The isomorphism was inferred from the classical isomorphism, the
agreement of Poisson brackets, and the $h$-adic completion.
In Section \ref{observable} we will show directly that the map $\Phi$ is 
an isomorphism.

\section{Temperley-Lieb theory}\label{temperley}
In this section we recall the correspondence between skeins and intertwiners
in the category of representations $\un{m}$ of $\mat$. 

Consider a rectangle $R=I\times I$ with $2n$ distinguished points; $n$
of them on $I\times \{0\}$ and $n$ on $I\times \{1\}$. Take the space
of blackboard framed  tangles with $n$ arc components ending at the
distinguished points. Its quotient by the Kauffman bracket skein
relations is denoted by $\kt(R,n)$. This quotient has an algebra
structure given by placing the bottom of a rectangle on the top of
another in such a way that the  distinguished points meet.

\begin{fact}\label{rect}
The algebra $\kt(R,n)$ is isomorphic to
$\mathrm{End}_{\mat}(\un{1}^{\otimes n})$, the space of  $\mat$-linear
maps of $\un{1}^{\otimes n}$ to itself. 
\end{fact} 
In the basis $\{e_{1/2}, e_{-1/2}\}$ of $\un{1}$,
the isomorphism is given by a tangle functor which makes the following
assignments. A local maximum is sent to the morphism 
$\mu :\un{1}\otimes \un{1} \rightarrow \un{0}$
 defined by 
\begin{equation}\label{max}
\mu (e_{1/2}\otimes e_{-1/2})=it, \quad \mu(e_{-1/2}\otimes
e_{1/2})=-it^{-1},    
\end{equation}
\[\mu(e_{-1/2}\otimes e_{-1/2})=\mu(e_{1/2}\otimes
e_{1/2})=0.\]
A local minimum is associated to the morphism
$\eta: \un{0} \rightarrow \un{1}\otimes \un{1}$ given by
\[\eta(1)=it e_{1/2}\otimes e_{-1/2}-it^{-1} e_{-1/2}\otimes
e_{1/2}. \]
\begin{fact}\label{joneswenzl}
The isomorphism takes the $n$-th Jones-Wenzl
idempotent to the intertwiner that  projects
$\un{1}^{\otimes n}$ onto its highest weight invariant subspace.
\end{fact}

Let $\mathbb{H}=\{(x,y)\ |\ y\geq 0\}$ be the closed upper half
plane. For any $n$, choose $2n$ distinguished points on the $x$-axis,
$\{(1,2,\dots,2n)\}$, 
and form a space of blackboard-framed tangles with $n$ arc components
ending  at the distinguished points. The quotient of this space by the
Kauffman bracket skein relations is denoted $\kt(\mathbb{H},2n)$.
\begin{fact}
$\kt(\mathbb{H},2n)\cong\mathrm{Inv}\left(((\un{1}\otimes\un{1})^{\otimes
    n})'\right)$ 
\end{fact}
\noindent The isomorphism is given by the same tangle functor as for
Fact \ref{rect}. 

An admissible triple $(m,n,p)$ determines a skein in
$\kt(\mathbb{H},m+n+p)$ consisting of a Kauffman triad with all three
legs attached to the $x$-axis. Fact \ref{joneswenzl} gives a canonical
inclusion of $\un{m}\otimes \un{n} \otimes \un{p}$ into 
$\un{1}^{\otimes m}\otimes\un{1}^{\otimes n}\otimes\un{1}^{\otimes p}
\cong (\un{1}\otimes\un{1}) ^{\otimes (m+n+p)/2}$.
\begin{fact}\label{triad}
 An admissible triple $(m,n,p)$ --- equivalently, a
Kauffman triad --- corresponds to a nonzero vector in the $1$-dimensional
space $\mathrm{Inv}((\un{m} \otimes \un{n} \otimes \un{p})')$.
\end{fact}

\section{Observables and the Kauffman bracket skein algebra}\label{observable}

Our goal is to prove a theorem analogous to  \cite[Theorem 10]{lattice}, but
replacing power series by complex numbers.

\begin{theorem}\label{iso}
 Let $\Gamma$ be a lattice and let $F$ be its envelope. Assume that
$ t \in \mathbb{C}\setminus\{0\}$ is not a root of unity. The algebra of
observables of lattice  
gauge field theory on $\Gamma$ based on $(\mat,\qg)$ is isomorphic
to $K_t(F)$.
\end{theorem}
\proof
From \cite{lattice} we have an algebra map from $\bC\cL$  to
$\cO(\Gamma)$ taking a link to the corresponding Wilson operator. As
in \cite[Theorem 10]{lattice} this map descends to 
$\Phi : \kt(F) \to \cO(\Gamma)$.

The following description of $\Phi$ is implicit in \cite{lattice}.
Since $\cO(\Gamma)$ is a homeomorphism invariant of $F$,
we can assume that the lattice $\Gamma$ comes
from giving orientation and ciliation to a 
trivalent spine $\gamma$ of $F$.
Let $L$ be an element of the  basis ${\mathcal B}_2$ of $\kt(F)$
(i.e., a  link with a simple diagram).
Choose an orientation of $L$. 
To compute the image of $L$ under $\Phi$,
first perform the
composition:
\begin{equation}\label{phi}\begin{CD}
\bigotimes_e\mat @>\Delta>>\bigotimes_e\mat^{\otimes  n(e)}
@>\otimes_e\rho_1^{\otimes n(e)}>>
\bigotimes_e(\un{1}^*\otimes\un{1})^{\otimes  n(e)}.
\end{CD}\end{equation}
Here $n(e)$ is the number of
strands of $L$ running along the edge $e$ of $\Gamma$, and $\Delta$
means comultiplying $n(e)-1$ times in the factor corresponding to
$e$. 

Second, in each factor where the corresponding segment of $L$
runs against the orientation of the edge of $\Gamma$,
apply $i D\otimes i  D^{-1}
:\un{1}^*\otimes\un{1}\rightarrow \un{1}\otimes\un{1}^*$. 
Here, the morphism $D:\un{1}^*\to\un{1}$ 
is defined by 
\begin{equation}\label{D}
D\left(e^{1/2}\right) = ite_{-1/2},\quad
D\left(e^{-1/2}\right) = -it^{-1}e_{1/2}.
\end{equation}

\begin{figure}
\lmk\makebox[53pt]{\tiny{$\phi(x)$}}
\makebox[53pt]{\tiny{$\phi(K^2x)$}}\rmk
\vspace{27pt}
\lmk\makebox[53pt]{\tiny{$\phi$}\hspace{10pt}\tiny{$x$}}
\makebox[53pt]{\tiny{$x$}\hspace{10pt}\tiny{$\phi$}}\rmk
\vspace{-36pt}
\lmk\epsfig{file=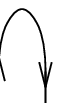}\hspace{40pt}
\epsfig{file=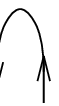}\rmk
\caption{Cap tangles}
\label{caps}
\end{figure}

Third, treat each ciliated vertex 
as a half plane with the cilium at infinity.
Up to isotopy, the link $L$ now appears as a collection of oriented
caps in each half plane.  
The cap pictured on the left of Figure
\ref{caps} is associated with the map $\un{1}^*\otimes\un{1}\to\bC$
where $\phi\otimes x\mapsto \phi(x)$, and $\un{1}\otimes\un{1}^*\to\bC$
is given by $x\otimes\phi\mapsto\phi(K^2x)$    for the cap on the right.

Finally, to obtain $\Phi(L)$, multiply the result by $(-1)^{|L|}$, 
where $|L|$ denotes the number of components of $L$.

Notice that $\rho_1$ sends the switch map of \cite{lattice} to the map
$i D\otimes i  D^{-1}$ and sends multiplications to
contractions. Hence our description of $\Phi(L)$ for $L\in{\mathcal B}_2$
 coincides with the Wilson operator. 
It follows from \cite{lattice} that $\Phi$ does not depend on the
choice of orientation of $L$. Extend it linearly to all of $\kt(F)$.
It is a homomorphism of algebras $\kt(F)$ and $\cO(\Gamma)$.
In order to prove that $\Phi$ is an isomorphism we factor it into two
maps which are isomorphisms on the level of vector spaces. 

The first map, expressed in the basis ${\mathcal B}_2$, is given
by a diagonal matrix with $1$'s and $-1$'s on the diagonal. 
The second map,
\[\Phi_u:\kt(F) \to \cO(\Gamma),\]
does not require a choice of orientation of a link $L$.
To compute the image of $L\in{\mathcal B}_2$ under $\Phi_u$ 
first perform  the composition  (\ref{phi}).
Second, apply the map $D$ to each copy of $\un{1}^*$.
Third, treat each vertex as a half plane and
associate the map $\mu$ from equation (\ref{max}) to
the (unoriented) caps.

Checking all possible orientations of $L$ and $\Gamma$ shows that 
$\Phi(L)=\pm\Phi_u(L)$.

By Fact \ref{triad}, the map $\Phi_u$ takes an
element of the basis 
${\mathcal B}_1$ (i.e., a skein obtained by an admissible coloring of 
$\gamma$)
to an observable coming from coloring the edges of $\Gamma$ with
corresponding irreducible representations of $\mat$. Thus, by 
Proposition \ref{basis}, the map $\Phi_u$
takes the basis ${\mathcal B}_1$ of $\kt(F)$ to a basis  of $\cO(\Gamma)$.

As $\Phi$ and $\Phi_u$ differ by a composition with an
isomorphism, both maps are isomorphisms.\qed

\end{document}